\documentclass[aap,MSNbibl,seceqn,citesort,dvips]{arximspdf}
\usepackage{mathrsfs}
\usepackage{graphicx}

\doi{10.1214/09-AAP664}
\volume{20}
\issue{4}
\pubyear{2010}
\firstpage{1205}
\lastpage{1218}

\makeatletter
\DeclareMathAlphabet\mathcaligr{OMS}{cmsy}{m}{n}
\renewcommand{\mathcal}{\mathcaligr}
\renewcommand{\emptyset}{\varnothing}
\newcommand{\eqref}[1]{(\ref{#1})}

\newtheorem{theorem}{Theorem}[section]
\newtheorem{proposition}{Proposition}[section]
\newtheorem{lemma}{Lemma}[section]

\newproclaim{remark}{Remark}[section]
\newproclaim{Question}{Question}
\makeatother

\begin{document}
\begin{frontmatter}

\title{Uniform convergence for complex $[\mathbf{0,1}]$-martingales}
\runtitle{Uniform convergence for complex $[0,1]$-martingales}

\begin{aug}
\author[A]{\fnms{Julien} \snm{Barral}\corref{}\ead[label=e1]{barral@math.univ-paris13.fr}},
\author[B]{\fnms{Xiong} \snm{Jin}\ead[label=e2]{Xiong.Jin@inria.fr}}
\and
\author[C]{\fnms{Beno\^{i}t} \snm{Mandelbrot}\ead[label=e3]{Benoit.Mandelbrot@yale.edu}}
\runauthor{J. Barral, X. Jin and B. Mandelbrot}
\affiliation{Universit\'e Paris 13, INRIA and Yale University}
\address[A]{J. Barral\\
LAGA (UMR 7539), Institut Galil\'ee\\
Universit\'{e} Paris 13\\
99 av. Jean-Baptiste Cl\'ement\\
93430 Villetaneuse\\
France \\
\printead{e1}} 
\address[B]{X. Jin\\
INRIA\\
Domaine de Voluceau\\
78153 Le Chesnay cedex\\
France\\
\printead{e2}}
\address[C]{B. Mandelbrot\\
Mathematics Department\\
Yale University\\
New Haven, Connecticut 06520\\
USA\\
\printead{e3}}
\end{aug}
\pdfauthor{Julien Barral, Xiong Jin, Benoit Mandelbrot}

\received{\smonth{2} \syear{2009}}
\revised{\smonth{11} \syear{2009}}

%
\begin{abstract}
Positive $T$-martingales were developed as a general framework that
extends the positive measure-valued martingales and are meant to
model intermittent turbulence. We extend their scope by allowing the
martingale to take complex values. We focus on martingales
constructed on the interval $T=[0,1]$ and replace random measures
by random functions. We specify a large class of such martingales
for which we provide a general sufficient
condition for almost sure uniform convergence to a nontrivial
limit. Such a limit yields new examples of naturally generated
multifractal processes that may be of use in multifractal signals
modeling.
\end{abstract}

%

\end{frontmatter}

\section{Introduction}\label{intro}

\subsection{Foreword about multifractal functions} Multifractal
analysis is a
natural framework to describe the
heterogeneity that is reflected in the distribution at small scales
of the H\"older singularities of a given locally bounded function or
signal $F \dvtx I\mapsto\mathbb{C}$ where $I$ is an interval. The H\"older
singularity of $F$ can be defined, at every point $t$, by
%
\[
h_F(t)=\liminf_{r\to0^+} \frac{\log\operatorname
{Osc}_F([t-r,t+r])}{\log(r)}
\]
or
\[
h_F(t)=\liminf_{n\to\infty
} \frac{\log_2\operatorname{Osc}_F(I_n(t))}{-n},
\]
where $I_n(t)$ is the dyadic interval of length $2^{-n}$ containing
$t$ and $\operatorname{Osc}_F(J)=\sup_{s,t\in J}|F(t)-F(s)|$. The
multifractal analysis of $F$ classifies points according to
$h_F$. It may compute the singularity spectrum
of $F$, that is, the Hausdorff dimension of the sets
$h_F^{-1}(\{h\})$
for $h\ge0$ or, more roughly, measure the asymptotic number of dyadic
intervals of generation $n$ needed to cover the sets $h_F^{-1}(\{h\})$ by
estimating the large deviation spectrum
%
%
\[
\label{large}
L_F(h)=\lim_{\varepsilon\to0}\limsup_{n\to\infty}\frac{\log_2\#
\{J\in\mathcal{G}_n, 2^{-n(h+\varepsilon)}\le\operatorname
{Osc}_F(J)\le2^{-n(h-\varepsilon)}\}}{n},
\]
where $\mathcal{G}_n$ is the set of dyadic intervals of generation $n$.
One says that $F$ is monofractal if there exists a unique $h\ge0$
such that $E_F(h)\neq\emptyset$ or $L_F(h)\neq- \infty$. Otherwise,
$F$ is multifractal
(see \cite{JAFFJMP,Riedi} for more details).

\subsection{Motivations and methods to build multifractal processes}
The main motivation for constructing and studying multifractal
functions or stochastic processes comes from the need to model
empirical signals for which the estimation of $L_{F}$ and related quantities
reveals striking scaling invariance properties. These signals concern
physical or social
intermittent phenomena like energy dissipation in turbulence
\cite{M1,M3,FrPa}, spatial rainfall \cite{Gupta}, human heart rate
\cite{Stanley}, internet traffic \cite{RRJLV} and
stock exchange prices \cite{Mandfin}. Models of these phenomena are
the statistically self-similar measures constructed in \cite
{M3,KP,BM2,BacryMuzy}. These objects are special examples of limit of
``$T$-martingales,'' which consist in a class of random measures developed in \cite{K2,K3}
after
the seminal work \cite{M1} about Gaussian multiplicative chaos
(see also \cite{Fan4,WaWi}). When $T=[0,1]$, these martingales and their
limit are also
used to build models of nonmonotonic scaling invariant signals as
follows: By
performing a multifractal time change in Fractional Brownian motions
or stable L\'evy processes
\cite{Mandfin,BacryMuzy,Riedi}, by integrating a
positive $[0,1]$-martingale with respect to the Brownian motion or
using such a martingale to specify the covariance of some Gaussian
processes to get new types of multifractal random walks
\cite{BacryMuzy,Ludena}, or by considering random wavelet series~whose coefficients are built from a multifractal measure \cite{ABM,BSw}.

\subsection{A natural alternative construction} This paper considers
the natural alternative to these constructions which allows the
multiplicative processes involved in
$[0,1]$-martingales to take complex values.

Let us now recall what are $[0,1]$-martingales.
Let $(\Omega,\mathcal{B},\mathbb{P})$ be a probability space, endow
the interval $[0,1]$ with the Borel $\sigma$-algebra
$\mathcal{B}([0,1])$ and the product space $[0,1]\times\Omega$ with the
product $\sigma$-algebra $ \mathcal{B}([0,1])\otimes\mathcal{B}$.
Let $(\mathcal{B}_n)_{n\ge1}$ be a nondecreasing sequence of
$\sigma$-algebras in $\mathcal{B}$. Also let $(Q_n)_{n\ge1}$ be a
sequence of complex-valued measurable functions defined on
$[0,1]\times\Omega$ such that for each $t\in[0,1]$,
$\{Q_n(t,\cdot),\mathcal{B}_n\}_{n\ge1}$ is a martingale of
expectation 1. Such a sequence of functions is called a
$[0,1]$-martingale. Given a Radon measure $\lambda$ on $[0,1]$, for
every $n\ge1$ we can define the random complex measure $\mu_n$
whose density with respect to $\lambda$ is equal to $Q_n$.

If the functions $Q_n$ take nonnegative values, then, with
probability 1, the sequence of Radon measures $(\mu_n)_{n\ge1}$
weakly converges to a measure $\mu$ (\cite{K2,K3}). This property is
an almost straightforward consequence of the positive martingale
convergence theorem and Riesz's representation theorem. When the
random functions $Q_n$ cease to be nonnegative, the martingales
$Q_n(t)$ need not be bounded in $L^1$ norm; hence the total
variations of the complex measures $\mu_n$ may diverge, and
$(\mu_n)_{n\ge1}$ need not converge almost surely weakly to an
element of the dual of $\mathcal{C}([0,1])$, the space of continuous
complex-valued functions over $[0,1]$.
In this paper we rather consider the sequence of random continuous
functions
\[
F_n \dvtx t\in[0,1]\mapsto\mu_n([0,t])=\int_0^tQ_n(u) \,\mathrm{d}\lambda(u).
\]
Then the following questions arise naturally:

\begin{Question} Does there exist a general necessary and sufficient
condition under which $(F_n)_{n\ge1}$
converges almost surely uniformly to a limit which is nontrivial
(i.e., different from 0) with positive probability?
\end{Question}
\begin{Question} When the sequence $(F_n)_{n\ge1}$ diverges, or
converges to 0 in $\mathcal{C}([0,1])$, can a natural normalization
of $F_n$ make it converge to a nontrivial multifractal limit
$\widetilde F$, at least in distribution?
\end{Question}
\begin{Question} Consider the case of strong or weak convergence to
a limit process $F$ or $\widetilde F$ having scaling invariance
properties. What is the multifractal nature of
$F$ (or $\widetilde F$), and does $F$ or $\widetilde F$ possess the
remarkable property to be naturally decomposed as a monofractal
function in multifractal time, like for some other classes of
multifractal functions \cite{Mandfin,JafMand2,BacryMuzy,Sadv}?
\end{Question}

\begin{figure}

\includegraphics{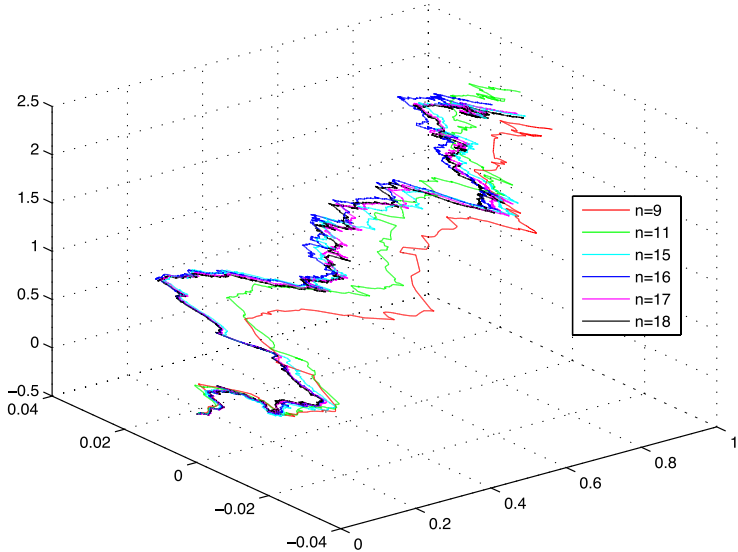}

\caption{A complex valued canonical dyadic cascade $F_{n}$ for $n=9$,
11, 15, 16, 17, 18.}\label{fig1}
\end{figure}

We will introduce a subclass of complex $[0,1]$-martingales,
namely $\mathcal{M}$, such that for $(Q_n)_{n\ge1}\in\mathcal{M}$,
we have a general sufficient condition for the almost sure uniform
convergence of $(F_n)_{n\ge1}$ to a nontrivial limit, as well as a
result of global H\"older regularity for the limit function
(Theorem~\ref{thm:1}). Our result makes it possible to
construct the complex extensions of some fundamental examples of
statistically self-similar positive multiplicative cascades
mentioned above (see
Section~\ref{ex} and an illustration in Figure~\ref{fig1}).

Companion papers \cite{BJ} and \cite{BJMpartII} provide further
results and answers to the previous questions in the particular case of
complex $b$-adic independent cascades (it is worth noting that these
objects also play a role in the study of directed polymers in a random
medium \cite{DES}).

Section~\ref{01mart} introduces the class $\mathcal{M}$, states
Theorem~\ref{thm:1} and provides fundamental examples in
$\mathcal{M}$. Section~\ref{proofth1} provides the proof of
Theorem~\ref{thm:1}. We end this section with some definitions.

\subsection{Definitions} Given an integer $b\geq2$, we denote by
$\mathscr{A}$ the alphabet $\{0,\ldots,b-1\}$ and define $\mathscr
{A}^*=\bigcup_{n\ge0} \mathscr{A}^n$ (by convention $\mathscr{A}^0$
is the set reduced to the empty word denoted $\emptyset$). For every
$n\ge0$, the length of an element of $\mathscr{A}^n$ is by definition
equal to $n$, and we denote it by $|w|$.
For $w \in\mathscr{A}^*$, we define $t_w=\sum_{i=1}^{|w|} w_i
b^{-i}$ and $I_w=[t_w,t_w+b^{-|w|}[$. For $n\ge1$ we define $T_n=\{
t_w: w\in\mathscr{A}^n\}\cup\{1\}$ and then $T_*=\bigcup_{n\geq1} T_n$.

For any $t\in[0,1)$ and $n\ge1$, we denote by $t|n$ the unique word
in $\mathscr{A}^n$ such that $t\in I_{t|n}$. We also denote by $t|0$
the empty word.

If $f\in\mathcal{C}([0,1])$ we denote by $\|f\|_\infty$ the norm
$\sup_{t\in[0,1]} |f(t)|$.

We denote by $(\Omega,\mathcal{B},\mathbb{P})$ the probability space
on which the random variables considered in this paper are defined.

\section{A class of complex $[0,1]$-martingales}\label{01mart}
\subsection{Definition} Consider a sequence of measurable complex functions
\[
P_n \dvtx \bigl([0,1]\times\Omega,\mathcal{B}([0,1])\otimes\mathcal{B}\bigr)
\mapsto(\mathbb{C},\mathcal{B}(\mathbb{C})),\qquad n\ge1.
\]

For $n\ge1$ and $I$, a subinterval of $[0,1]$, let $\mathcal{F}_n^I$
be the $\sigma$-field generated in $\mathcal{B}$ by the family of
random variables $\{P_{m}(t,\cdot)\}_{t\in I, 1\le m\le n}$. Also let
$\overline{\mathcal{F}}{}^{I}_{n}$ be the $\sigma$-field generated in
$\mathcal{B}$ by the family of random variables $\{P_m(t,\cdot)\}_{t
\in I,m>n}$. The $\sigma$-fields $\mathcal{F}^{[0,1]}_{n}$ and
$\overline{\mathcal{F}}{}^{[0,1]}_{n}$ are simply denoted by $\mathcal
{F}_{n}$ and $\overline{\mathcal{F}}_{n}$.

\begin{longlist}[(P3)]
\item[(P1)] For all $t\in[0,1]$, $P_n(t,\cdot)$ is integrable, and
$\mathbb{E}(P_n(t,\cdot))=1$.
\item[(P2)] For every $n\ge1$, $\mathcal{F}_n$ and $\overline
{\mathcal{F}}_n$ are independent.
\item[(P3)] There exist two integers $b\ge2$ and $N\geq1$ such
that for every $n\ge1$ and every family $\mathcal{G}$ of $b$-adic
subintervals of $[0,1]$ of generation $n$ such that ${d}(I,J) \geq N
b^{-n}$ for every $I \neq J \in\mathcal{G}$, the $\sigma$-algebra's
$\overline{\mathcal{F}}{}^{I}_n$, $I\in\mathcal{G}$, are mutually
independent, where $d(I,J)=\inf\{|t-s| \dvtx s\in I, t\in J\}$.
\end{longlist}

Under the properties (P1) and (P2), for each $t\in
(0,1)$ the sequence
\[
Q_n(t,\cdot)=\prod_{k=1}^nP_k(t,\cdot)
\]
is a martingale of expectation 1 with respect to the filtration $\{
\mathcal{F}_n\}_{n\ge1}$.

We denote by $\mathcal{M}$ the class of martingales $(Q_n)_{n\ge1}$
obtained as above and which satisfy properties (P1)--(P3).

We denote by $\mathcal{M'}$ the subclass of $\mathcal{M}$ of those
$(Q_n)_{n\ge1}$ which, in addition to (P1)--(P3),
satisfy the statistical self-similarity property:

\begin{longlist}
\item[(P4)] Let $b$ be as in (P3). For every closed
$b$-adic subinterval $I$ of $[0,1]$, let $n(I)$ and $S_I$,
respectively, stand for the generation of $I$ and the canonical affine
map from $[0,1]$ onto $I$. The processes $ (P_{n(I)+n}\circ S_I )_{n\ge
1}$ and $(P_{n})_{n\ge1}$ have the same distributions.
\end{longlist}

Let $\lambda$ be a Radon measure on $[0,1]$.
If $(Q_n)_{n\ge1}\in\mathcal{M}$, for $n\ge1$, we define
%
\begin{equation}\label{Fn}
F_n(t)=\int_0^t Q_n(u) \,\mathrm{d}\lambda(u).
\end{equation}

\subsection{Convergence theorem for $(F_n)_{n\ge1}$}

Theorem~\ref
{thm:1} provides a sufficient condition for the almost sure uniform
convergence of $F_n$, as $n$ tends to $\infty$, to a limit $F$ such
that $\mathbb{P}(F\neq0)>0$. This condition is the extension of the
condition introduced in Part II of \cite{BMS2} to show that when
$(Q_n)_{n\ge1}$ is nonnegative, the sequence of measures $F_n'$
converge almost surely weakly to a random measure $\mu$ such that
$\mathbb{P}(\mu\neq0)>0$. When $(Q_n)_{n\ge1}$ is not nonnegative,
the uniform convergence of $F_n$ is a more delicate issue.

For $p\in\mathbb{R}_+$ and $n\ge1$ we define
%
\begin{eqnarray}\label{eqn:1}
S(n,p)&=&\sum_{w\in\mathscr{A}^n} \lambda(I_{w})^{p-1}\int_{I_{w}}
\mathbb{E}(|Q_n(t)|^p) \,\mathrm{d}\lambda(t)\quad\mbox{and}
\\\label{eqn:2}
\varphi(p)&=&\liminf_{n\rightarrow\infty} \biggl(-\frac{1}{n} \log_b S(n,p)\biggr).
\end{eqnarray}
We notice that $\varphi$ is a concave function of $p$, $\varphi(0)\le
0$ by construction, and that due to our assumption that $\mathbb
{E}(Q_n(t))=1$, we also have $\varphi(1)\le0$.

\eject
\begin{theorem}\label{thm:1}
\begin{longlist}[(2)]
\item[(1)] Suppose that $\varphi(p)>0$ for some $p\in(0,1)$, and that
there exists a function $\psi \dvtx\mathbb{N}_+\to\mathbb{R}_+$ such
that $\psi(n)=o(n)$ and $\mathbb{E}(\sup_{t\in I_w}|Q_n(t)|^p)\le
\exp(\psi(n))\times\break \mathbb{E}(|Q_n(t)|^p)$ for all $n\ge1$, $w\in
\mathscr{A}^n$ and $t\in I_w$. Then, with probability 1, $F_{n}$
converge uniformly to 0 as $n\to\infty$.
\item[(2)]
Let $p\in(1,2]$. Suppose that $\varphi(p)>0$. The functions $F_{n}$
converge uniformly, almost surely and in $L^1$ norm, to a limit $F$, as
$n\to\infty$. The function $F$ is $\gamma$-H\"older continuous for
all $\gamma\in(0,\max_{q\in(1,p]}\varphi(q)/q)$. Moreover,
$\mathbb{E}(\|F\|_\infty^p)<\infty$.
\end{longlist}
\end{theorem}

\begin{remark}
%
(1) The proof of Theorem \ref{thm:1}(1) will show that this result
does not require \textup{(P1)},
\textup{(P2)} or \textup{(P3)}.
The existence of the function $\psi$ corresponds to a kind of bounded
distortion principle.
%


(2) Under the assumptions of Theorem~\ref{thm:1}(2), let $\beta=\min
\{p\in[1,2) \dvtx\varphi(p)=0\}$. The nonnegative sequence $(Q^{(\beta
)}_n)_{n\ge1}=(|Q_n|^\beta)_{n\ge1}$ is an element of $\mathcal
{M}$, and by construction, the corresponding function $\varphi$ is
positive near $1^+$. Consequently, the sequence $F_n^{(\beta)}$
defined by $\int_0^\cdot Q^{(\beta)}_n(u) \,\mathrm{d}u$ converges
uniformly to a nondecreasing function $F^{(\beta)}$. Inspired by the
results obtained in \cite{BJMpartII}, it is natural to ask under which
additional assumptions it is possible to write $F=B_{1/\beta}\circ
F^{(\beta)}$ where $B_{1/\beta}$ is a monofractal function of
exponent $1/\beta$.

(3) Suppose that $\varphi$ is not positive over $[0,2]$. In the case
where the martingale $(F_n(1))_{n\ge1}$ is not bounded in $L^2$ norm,
inspired again by what is done in \cite{BJMpartII}, it is natural to
look at the process $F_n/\sqrt{\mathbb{E}(F_n(1)^2)}$ and seek for
conditions under which it converges in distribution, as $n\to\infty$.
\end{remark}

\subsection{Examples}\label{ex}
\subsubsection*{Homogeneous $b$-adic independent cascades}
We consider the complex extension of the nonnegative
$[0,1]$-martingales introduced in \cite{M3}. Let $b$ be an integer
$\ge2$ and for every $k\ge0$ let $W^{(k)}=(W^{(k)}_0,\ldots
,W^{(k)}_{b-1})$ be a vector such that each~of its components is
complex, integrable and has an expectation equal to 1. Then, consider
$\{W^{(|w|)}(w)\}_{w\in\mathscr{A}^*}$, a family of independent
vectors such that for each $k\ge0$ and $w\in\Sigma_k$ the vector
$W^{(k)}(w)$ is a copy of $W^{(k)}$.

An element of $\mathcal{M}$ is obtained as follows. For $t\in[0,1)$
and $n\ge1$ let
$
P_n(t)=W^{(n-1)}_{t_{n}}(t|n-1) \mbox{ and then }Q_n(t)=\prod
_{k=1}^nP_k(t)$. If $\lambda$ is the inhomogeneous Bernoulli measure
associated with a sequence of probability vectors $(\lambda
^{(k)}=\lambda^{(k)}_0,\ldots,\lambda^{(k)}_{b-1})_{k\ge0}$, then
\[
\varphi(p)=\liminf_{n\to\infty}\Biggl(-\frac{1}{n}\sum_{k=0}^{n-1}\log
_b\mathbb{E} \Biggl(\sum_{i=0}^{b-1} \bigl(\lambda^{(k)}_i \big|W^{(k)}_i\big| \bigr)^p \Biggr)\Biggr).
\]
If all the vectors $W^{(k)}$ have the same distribution as a vector
$W$, then $(Q_n)_{n\ge1}$ belongs to $\mathcal{M}'$.
\textit{Canonical cascades} correspond to $W$ whose components are
i.i.d. and $\lambda$ equal to the Lebesgue measure. Then a necessary
and sufficient condition for the almost sure uniform convergence of
$F_n$ to a nontrivial limit is $\varphi'(1^-)>0$ if $W\ge0$ \cite
{KP,KAIHP} and $\varphi(p)>0$ for some $p\in(1,2]$ for the special
``monofractal'' examples considered in \cite{BM}.


\subsubsection*{Compound Poisson cascades} We consider the complex
extension of the nonnegative $[0,1]$-martingales introduced in \cite
{BM2}. Let $\nu$ be a positive Radon measure over $(0,1]$ and denote
by $\Lambda$ the measure $\mbox{Leb}\otimes\nu$ where $\mbox{Leb}$
stands for the Lebesgue measure over $\mathbb{R}$. We consider a
Poisson point process $S$ of intensity $\Lambda$. To each point $M$ of
$S$, we associate a random variable $W_M$ picked in a collection of
random variables that are independent, independent of $S$, and are
identically distributed with an integrable complex random variable $W$.
We fix $\beta>0$, and for $n\ge1$ and $t\in[0,1]$ we define the
truncated cone
\[
\Delta\mathcal{C}_n (t)= \{(t',r) \dvtx b^{-n}<r\le b^{1-n}, t-\beta
r/2\le t'<t+\beta r/2 \}.
\]
We obtain an element of $\mathcal{M}$ as follows. For $t\in[0,1)$ and
$n\ge1$ we define
\[
P_n(t)= e^{-\Lambda (\Delta\mathcal{C}_n(t) )(\mathbb{E}(W)-1)
}\prod_{M\in S\cap\Delta\mathcal{C}_n(t)}W_M,
\]
 and then
 \[
 Q_n(t)=\prod_{k=1}^nP_k(t).
\]
If $\lambda$ is the Lebesgue measure and $\widetilde\beta$ stands
for $\limsup_{n\to\infty}n^{-1} \log_b \Lambda(\bigcup
_{k=1}^n\Delta\mathcal{C}_{n})$,
\[
\varphi(p)=p -1+\widetilde\beta \bigl(p\bigl(\mathbb{E}(\Re W)-1\bigr)-\bigl(\mathbb
{E}(|W|^p)-1\bigr) \bigr).
\]
If, moreover, there exists $\delta>0$ such that $\nu
(\mathrm{d}r)=\delta\,\mathrm{d}r/r^2$, that is, if $\Lambda$ possesses
scaling invariance properties, we have $\widetilde\beta=\beta\delta
$, and $(Q_n)_{n\ge1}$ belongs to $\mathcal{M}'$.



\subsubsection*{Log-infinitely divisible cascades} This example is an
extension of compound Poisson cascades when the weights $W_M$ take the
form $\exp(L_M)$, and, in particular, the $W_M$ do not vanish.
We use the notations of the previous section and take $\beta=\delta
=1$. Let $\psi$ be a characteristic L\'evy exponent $\psi$ defined on
$\mathbb{R}^2$, that is,
%
\begin{equation}\label{psi}
\qquad\psi \dvtx\xi\in\mathbb{R}^2\mapsto i\langle\xi|a \rangle-Q(\xi
)/2+\int_{\mathbb{R}^2} \bigl(1-e^{i\langle\xi|x\rangle}+i\langle\xi
|x\rangle\mathbf{1}_{|x|\le1}\bigr)  \pi(\mathrm{d}x),
\end{equation}
where
$a\in\mathbb{R}^2$, $Q$ is a nonnegative quadratic form and $\pi$
is a Radon measure on $\mathbb{R}^2\setminus\{0\}$ such that $\int
(1\land|x|^2)\pi(\mathrm{d}x)<\infty$.

Then let $\rho=(\rho_1,\rho_2)$ be an independently scattered
infinitely divisible random $\mathbb{R}^2$-valued measure on $\mathbb
{R}\times\mathbb{R}^*_+$ with $\Lambda$ as control measure and $\psi
$ as L\'evy exponent (see \cite{RR} for the definition). In
particular, for every Borel set $B\in\mathbb{R}\times\mathbb
{R}^*_+$ and $\xi\in\mathbb{R}^2$ we have
\[
\mathbb{E}\bigl(e^{i\langle\xi|\rho(B)\rangle}\bigr)=\exp (\psi(\xi
)\Lambda(B) ),
\]
and for every finite family $\{B_i\}$ of pairwise disjoint Borel
subsets of $ \mathbb{R}\times\mathbb{R}^*_+$ such that $\Lambda
(B_i)<\infty$, the random variables $\rho(B_i)$ are independent.

Let $I_1$ be the interval of those $\xi_1\in\mathbb{R}$ such that
$\int_{|x|\ge1} e^{ \xi_1x_1}\pi(\mathrm{d}x)<\infty$. The
function $\psi$ has a natural extension $\widetilde\psi$ to
$\mathcal{D}=\mathbb{R}^2\cup(-i I_1\times\mathbb{R})$ given by
the same expression as in~(\ref{psi}) if we extend $Q$ to an Hermitian
form on $\mathbb{C}^2$. Then for every $\xi\in\mathcal{D}$ and
every Borel subset of $ \mathbb{R}\times\mathbb{R}^*_+$ we have
$\mathbb{E}(e^{i\langle\xi|\rho(B)\rangle})=\exp (\widetilde\psi
(\xi)\Lambda(B) )$.

Now, we assume that $\xi_0=(-i,1)\in\mathcal{D}$, and without loss
of generality we set
\[
\widetilde\psi:=\widetilde\psi- \widetilde\psi(\xi_0).
\]
Then, with the same definition of cones as in the previous section, if
$n\ge1$ and $t\in[0,1]$, we define
\[
P_n(t) =\exp [ \langle\xi_0|\rho (\Delta\mathcal{C}_n(t) )
\rangle ]
=\exp [\rho_1 (\Delta\mathcal{C}_n(t) )+i\rho_2 (\Delta\mathcal
{C}_n(t) ) ]
\]
and $ Q_n(t)=\prod_{k=1}^nP_n(t)$. If we take $\lambda$
equal to the Lebesgue measure, and if $p\in\mathbb{R}$ is such that
$(-ip,0)\in\mathcal{D}$, then
%
\begin{equation}\label{phiidc}
\varphi(p)=p-1-\widetilde\beta\widetilde\psi(-ip,0).
\end{equation}

In the positive case, this construction that has been proposed has an
extension of compound Poisson cascades in \cite{BacryMuzy}. If $\nu
(\mathrm{d}r)=\mathrm{d}r/r^2$, then $(Q_n)_{n\ge1}$
belongs to $\mathcal{M}'$. In \cite{BacryMuzy}, a modification of
$P_1(t)$ is introduced, which yields a nice exact statistical scaling
invariance property for the increments of the limit measure. It can be
easily checked that this property, which is different from the
statistical self-similarity imposed by~(P4), also holds for the complex
extension.


%

%


%

%


\section[Proof of Theorem 2.1]{Proof of Theorem~\protect\ref{thm:1}}\label{proofth1}\leavevmode

\begin{pf*}{Proof of Theorem~\ref{thm:1}(1)}
For any $w\in\mathscr{A}^*$ and $n\ge1$, define
%
\begin{equation}
\Delta F_n(I_w) =F_n(t_w+b^{-n})-F_n(t_w)=\int_{I_{w}}Q_n(t) \,\mathrm{d}\lambda(t).
\label{eqn:6}
\end{equation}
We have
$
\mathbb{E}(\|F_n\|_\infty^p)\le\mathbb{E} ( (\sum_{w\in\mathscr
{A}^n} |\Delta F_n(I_w)| )^p )
\le\mathbb{E} (\sum_{w\in\mathscr{A}^n} |\Delta F_n(I_w)|^p ),
$\break
where we have used the subadditivity of $x\ge0\mapsto x^p$ ($p\in
(0,1]$). Thus
\begin{eqnarray*}
\mathbb{E}(\|F_n\|_\infty^p)&\le&\sum_{w\in\mathscr{A}^n} \mathbb{E}
\biggl( \bigg|\int_{I_w}Q_n(t) \,\mathrm{d}\lambda(t) \bigg|^p \biggr)\\
&\le&
\sum_{w\in\mathscr{A}^n} \lambda(I_w)^p \mathbb{E}\Bigl(\sup_{t\in
I_w}|Q_n(t)|^p\Bigr)\\
&\le& \sum_{w\in\mathscr{A}^n}\exp(\psi(n)) \lambda(I_w)^{p-1}
\int_{I_w}\mathbb{E}(|Q_n(t)|^p) \,\mathrm{d}\lambda(t)\\
&=&\exp(\psi(n)) S(n,p).
\end{eqnarray*}
Due to the property of $\psi(n)$, we have $\limsup_{n\to\infty}
\log_b (\mathbb{E}(\|F_n\|_\infty^p) )/n\le\break -\varphi(p)<0$. This
implies the result.
\end{pf*}

\begin{pf*}{Proof of Theorem~\ref{thm:1}(2)} The two following crucial
statements, which take natural and classical forms, will be proved at
the end of the section.

\begin{proposition}\label{controlbadic}
There exists a constant $C_p>0$ such that
%
\begin{equation}
(\forall n\ge2)\qquad \mathbb{E}\Bigl(\max_{t\in T_n}|F_n(t)-F_{n-1}(t)|^p\Bigr)\leq
C_p S(n,p).
\label{eqn:3}
\end{equation}
Consequently, for every $b$-adic number $t\in T_*$, $F_n(t)$ converges
almost surely and in $L^p$ norm as $n\to\infty$.
\end{proposition}

\begin{proposition}\label{holder}
Let $\gamma\in(0,\max_{q\in(1,p]}\varphi(q)/q)$. With probability
1, there exists $\eta_\gamma>0$ such that for any $t,s\in T_*$ such
that $|t-s|<\eta_\gamma$ we have
%
\begin{equation}
\sup_{n\ge1}|F_n(t)-F_n(s)| \leq C_\gamma|t-s|^\gamma,
\label{eqn:4}
\end{equation}
where $C_\gamma$ is a constant depending on $\gamma$ only.
\end{proposition}

Since $F_n(0)=0$ almost surely for all $n\ge1$, it follows from
Propositions~\ref{holder} and Ascoli--Arzela's theorem that, with
probability 1, the sequence of continuous functions $(F_n)_{n\ge1}$ is
relatively compact, and all the limits of subsequences of $F_n$ are
$\gamma$-H\"older continuous for all $0<\gamma<\max_{q\in
(1,p]}\varphi(q)/q$. Moreover, Proposition~\ref{controlbadic} tells
us that, with probability 1, $F_n$ is convergent over the dense
countable subset $T_*$ of $[0,1]$. This yields the uniform convergence
of $F_n$ and the H\"older regularity of the limit $F$.

We then prove that $\|\|F(t)\|_\infty\|_p<\infty$. For $n\ge1$, let
$M_{n}=\max_{t\in T_n} |F_{n}(t)|$. We have
%
\begin{equation}
M_{n+1}\leq M_n+\max_{t\in T^{n}} |F_{n+1}(t)-F_n(t)|+b\cdot\max
_{w\in\mathscr{A}^{n+1}} |\Delta F_{n+1}(w)|.
\label{eqn:11}
\end{equation}
Then Minkowski's inequality yields
\[
\|M_{n+1}\|_p\leq\|M_{n}\|_p+\Big\|\max_{t\in T_n} |F_{n+1}(t)-F_n(t)|\Big\|
_p+b\cdot\Big\|\max_{w\in\mathscr{A}^{n+1}} |\Delta F_{n+1}(w)|\Big\|_p.
\]
Also, due to Proposition~\ref{controlbadic} we have $\sum_{n\geq1} \|
\max_{t\in T_n} |F_{n+1}(t)-F_n(t)|\|_p<\infty.$ Moreover,
\[
\Big\|\max_{w\in\mathscr{A}^{n+1}} |\Delta F_{n+1}(w)|\Big\|_p\leq \biggl(\sum
_{w\in\mathscr{A}^{n+1}} \mathbb{E}(|\Delta F_{n+1}(w)|^p)
\biggr)^{1/p}\leq S(n+1,p)^{1/p},
\]
so $\sum_{n\geq1}\|\max_{w\in\mathscr{A}^{n+1}} |\Delta
F_{n+1}(w)|\|_p<\infty$. This implies $\sup_{n\ge1} \|M_{n}\|
_p<\infty$, and since, with probability 1, $F_n$ converges uniformly
to $F_\infty$, and $T_*$ is dense in $[0,1]$, we get $
\|\sup_{t\in[0,1]} |F(t)|\|_p\le\liminf_{n\rightarrow\infty} \|
M_{n}\|_p<\infty$. In particular, $F$~belongs to $L^1$ and for every
$n\ge1$, the conditional expectation of $F$ with respect to $\mathcal
{F}_n$ is well defined, and it converges almost surely and in $L^1$
norm to $F$ (see Proposition V-2-6 in \cite{Neveu}). It remains to
prove that $F_n=\mathbb{E}(F|\mathcal{F}_n)$ almost surely. For every
$t\in T_*$, we have shown that the martingale $(F_n(t),\mathcal
{F}_n)_{n\ge1}$ is uniformly integrable, so $F_n(t)=\mathbb
{E}(F(t)|\mathcal{F}_n)$ almost surely. Consequently, since $T_*$ is
countable, with probability 1, the restriction of $\mathbb
{E}(F|\mathcal{F}_n)$ coincides with the function $F_n$ over $T_*$.
Moreover, these two random functions are continuous and $T_*$ is dense
in $[0,1]$, so, with probability 1, they are equal.

\begin{pf*}{Proof of Proposition~\ref{controlbadic}}
Fix $n\geq2$ and
denote the elements of~$T_n$ by~$t_j$, $0\leq j \leq b^n$, where
$0=t_0<t_1<\cdots<t_{b^{n}}=1$. Also define $J_j=[t_j,t_{j+1}]$ for
$0\leq j <b^n$. We can write
\[
F_{n}(t_j)-F_{n-1}(t_j)=\sum_{k=0}^{j-1} \int_{J_k} U(t)V(t) \,\mathrm{d}\lambda(t)
\]
with $U(t)=Q_{n-1}(t)$ and $V(t)=P_n(t)-1$. Then we divide the family
$\{J_j\}_{0\leq j< b^n}$ into $bN$ sub-families, namely the $\{J_{bNk+i
}\}_{k\geq0 , 0\leq bNk+i <b^n}$, for $0\leq i \leq bN-1$. Also we
define $M_n=\max_{0\leq j\le b^n} |F_{n}(t_j)-F_{n-1}(t_j)|$ and
remark that
\begin{eqnarray*}
M_n\leq bN \mathop{\max_{0\leq j<b^n}}_{ 0 \leq i \leq bN-1} \bigg|\mathop
{\sum_{k\geq0}}_{0\leq bNk+i \leq j} \int_{J_{bNk+i}} U(t)V(t) \,\mathrm{d}\lambda(t) \bigg|.
\end{eqnarray*}
By raising both sides of the previous inequality to the power $p$ we
can get
%
\begin{eqnarray}\label{previneq}
 M_n^p
&\leq& (bN)^p \mathop{\max_{0\leq j<b^n}}_{0 \leq i
\leq bN-1} \bigg|\int_{J_{bNk+i}} U(t)V(t) \,\mathrm{d}\lambda(t) \bigg|^p\nonumber\\[-8pt]\\[-8pt]
&\le&(bN)^p\sum_{i=0}^{bN-1} \max_{0\leq j\leq b^n}
\bigg|\mathop{\sum_{k\geq0}}_{ 0\leq bNk+i \leq
j} \int_{J_{bNk+i}} U(t)V(t) \,\mathrm{d}\lambda(t) \bigg|^p.\nonumber
\end{eqnarray}
We are going to use the following lemma. It is proved for real valued
random variables in \cite{Shao1}, and its extension to the complex
case is immediate.

\begin{lemma}\label{lemma:2}
Let $p\in(1,2]$. There exists a constant $C_p>0$ such that for every
$n\ge1$ and every sequence $\{V_j\}_{1\leq j \leq n} $ of independent
and centered complex random variables we have
\[
\mathbb{E} \Biggl(\max_{1\leq k\leq n} \Bigg|\sum_{j=1}^{k}V_j\Bigg |^p \Biggr)\leq C_p
\sum_{j=1}^{n} \mathbb{E}(|V_j|^p).
\]
\end{lemma}

Due to (P3), for each $0\leq i \leq bN-1$, the restrictions
of the function $V(t)$
to the intervals $J_{bNk+i}$, $0\leq bNk+i<b^n$, are centered and independent.
Also, due to~(P2), the functions $U(t)$ and $V(t)$ are
independent. Consequently,
by taking the conditional expectation with respect to $\mathcal
{F}_{n-1}$ in~\eqref{previneq}
and using Lemma~\ref{lemma:2} we get for each $0 \leq i \leq bN-1$
\begin{eqnarray*}
&&\mathbb{E} \biggl(\max_{0\leq j\leq b^n} \bigg|\mathop{\sum_{k\geq0}}_{0\leq
bNk+i \leq j} \int_{J_{bNk+i}} U(t)V(t) \,\mathrm{d}\lambda(t) \bigg|^p
\Big|{\mathcal{F}}_{n-1} \biggr)\\
&&\qquad\le C_p
\mathop{\sum_{k\geq0}}_{ 0\leq bNk+i <b^n} \mathbb{E} \biggl( \bigg|\int
_{J_{bNk+i}} U(t)V(t) \,\mathrm{d}\lambda(t)\bigg |^p \Big|{\mathcal{F}}_{n-1} \biggr).
\end{eqnarray*}
This implies
%
\begin{equation}\label{Lp-conv}
\mathbb{E} (M_n^p |{\mathcal{F}}_{n-1} )\leq \widetilde C_p \sum
_{0\leq j\leq b^n} \mathbb{E} \biggl( \bigg|\int_{J_{j}} U(t)V(t) \,\mathrm{d}\lambda(t)\bigg |^p \Big|{\mathcal{F}}_{n-1} \biggr)
\end{equation}
with $\widetilde C_p=C_p (bN)^{p+1} $. Now, since $p>1$, the Jensen
inequality yields
\[
\bigg|\int_{I_{j}} U(t)V(t) \,\mathrm{d}\lambda(t) \bigg|^p \leq\lambda
(I_{j})^{p-1} \int_{I_{j}} |U(t)V(t)|^p \,\mathrm{d}\lambda(t).
\]
Moreover, since $\mathbb{E}(|P_{n}(t)|)\ge1$ and $p\ge1$, we have
%
\begin{equation}\label{momentp}
\mathbb{E}(|V(t)|^p)\leq2^{p-1}\bigl(1+\mathbb{E}(|P_{n}(t)|^p)\bigr)\leq2^p
\mathbb{E}(|P_{n}(t)|^p).
\end{equation}
Thus, taking the expectation in ($\ref{Lp-conv}$) yields
\[
\mathbb{E}\Bigl(\max_{t\in T_n}|F_n(t)-F_{n-1}(t)|^p\Bigr)\leq2^p\widetilde
C_p S(n,p),
\]
that is, (\ref{eqn:3}). If $\varphi(p)>0$, by definition of $\varphi
$, $S(n,p)$ converges exponentially fast to 0; hence the series $\sum
_{n\ge1}S(n,p)^{1/p}$ converge and, due to (\ref{eqn:3}) and the fact
that $T_*=\bigcup_{n\ge0}T_n$, $F_{n}(t)$ converges almost surely and
in $L^p$ norm as $n\rightarrow\infty$ for all $t\in T_*$.
\end{pf*}

\begin{pf*}{Proof of Proposition~\ref{holder}} Recall \eqref{eqn:6}.
Let $q\in(1,p]$ such that \mbox{$\varphi(q)>0$.} It follows from~(P1) that $(\Delta F_n(I_w) )_{n\ge1}$ is a martingale, so Doob's
and then Jensen's inequalities yield a constant $C_{q} $ such that for
$n\ge1$
\[
\mathbb{E}\Bigl(\max_{1\le k\le n}|\Delta F_k(I_w) |^{q}\Bigr)\leq C_{q}
\mathbb{E}(|\Delta F_n(I_w) |^{q})\le C_q \lambda(I_w)^{q-1} \int
_{I_{w}}\mathbb{E}(|Q_n(t)|^{q}) \,\mathrm{d}\lambda(t).
\label{eqn:7}
\]
Consequently
%
\begin{equation}
\sum_{w\in A^n} \mathbb{E}\Bigl(\max_{1\le k\le n}|\Delta F_k(I_w) |^{q}\Bigr)
\leq C_{q} S(n,{q}).
\label{eqn:8}
\end{equation}

By using Markov's inequality as well as (\ref{eqn:8}) and Proposition~\ref{controlbadic}, we get
\begin{eqnarray*}
&&
\mathbb{P}\Bigl(\max_{w\in A^n}\max_{0\leq k \leq n} |\Delta F_k(I_w) |
> b^{-n\gamma} \mbox{ or } \max_{t\in T_n} |F_{n}(t)-F_{n-1}(t)| >
b^{-n\gamma}\Bigr) \\
&&\qquad\leq \sum_{w\in A^n} \mathbb{P}\Bigl(\max_{0\leq k \leq n} |\Delta
F_k(I_w) | >b^{-n\gamma}\Bigr) + \mathbb{P}\Bigl(\max_{t\in T_n}
|F_{n}(t)-F_{n-1}(t)| > b^{-n\gamma}\Bigr)\\
&& \qquad\leq\sum_{w\in A^n} b^{n\gamma q} \cdot\mathbb{E}\Bigl(\max_{0\leq k
\leq n} |\Delta F_k(I_w) |^{q}\Bigr) + b^{n\gamma q} \cdot\mathbb{E}\Bigl(\max
_{t\in T_n} |F_{n}(t)-F_{n-1}(t)| ^{q }\Bigr)\\
&&\qquad\leq
C_{q} b^{n\gamma q} S(n,q),
\end{eqnarray*}
where $C_{q} $ is another constant depending only on $q$. Since $\gamma
\in(0,\varphi(q)/q)$, by definition of $\varphi(q)$ the series $\sum
_{n\ge1}b^{n\gamma q} S(n,q)$ converges, and by the Borel--Cantelli
lemma, with probability 1, there exists $n_1$ such that for all $n\geq n_1$,
%
\begin{equation}
\qquad\max_{w\in A^n}\max_{0\leq k \leq n} |\Delta F_k(I_w) |
\le
b^{-n\gamma} \quad\mbox{and}\quad \max_{t\in T_n} |F_{n}(t)-F_{n-1}(t)| \leq
b^{-n\gamma}.
\label{eqn:9}
\end{equation}

Now fix $n \geq n_1$. We are going to prove by induction that for all
$M\geq n+1$ and $t,s \in T_M $ such that $0<t-s<b^{-n}$ we have
%
\begin{eqnarray}\label{eqn:10}
\Delta_M(t,s) &\leq&2b \sum_{m=n+1}^{M} b^{-m\gamma}\qquad\mbox{where}\nonumber\\[-8pt]\\[-8pt]
\Delta_M(t,s) &=&\max_{0\leq k\leq M}|F_k(t)-F_k(s)|.\nonumber
\end{eqnarray}
If $M=n+1$, then there exist $i$ and $i'$ with $0<i-i'<2b$ such that
$t=ib^{-(n+1)}$ and $s=i'b^{-(n+1)}$, so due to~(\ref{eqn:9}) applied
at generation $n+1$,
\[
\Delta_{n+1}(t,s)\leq(i-i')b^{-(n+1)\gamma}\leq2b\cdot
b^{-(n+1)\gamma}.
\]
Now let $M\ge n+1$ and suppose that (\ref{eqn:10}) holds for all
$n+1\leq m \leq M$. Let $t,s \in T_{M+1}$ such that $0 < t-s< b^{-n}$.
If there is no element of $T_M$ between~$s$ and $t$, then (\ref
{eqn:9}) yields $\Delta_{M+1}(t,s)\leq(b-1)b^{-(M+1)\gamma}$.
Otherwise, consider $\bar t=\max\{u\in T_{M}:u\leq t\}$ and $
\bar s=\min\{u\in T_{M} \dvtx u\geq s\}$.
We have
\begin{eqnarray*}
s&\leq&\bar s\leq\bar t\leq t ,\qquad t-\bar t\leq(b-1) b^{-(M+1)} ,\\
 \bar
s-s&\leq&(b-1)b^{-(M+1)} ,\qquad \bar t-\bar s< b^{-n}.
\end{eqnarray*}
Since $\bar s$ and $\bar t$ belong to $T_{M}\subset T_{M+1}$, we deduce
from (\ref{eqn:9}) that
\[
\cases{
\max \{\Delta_{M+1}(t,\bar t), \Delta_{M+1}(\bar s,s)
\}\leq(b-1)b^{-(M+1)\gamma},\cr
\max \{|F_{M+1}(\bar s)-F_{M}(\bar s)|, |F_{M+1}(\bar
t)-F_{M}(\bar t)| \}\le b^{-(M+1)\gamma}.
}
\]
Also, due to (\ref{eqn:10}) we have $\Delta_{M}(\bar
t,\bar s) \leq2b \sum_{m=n+1}^{M} b^{-m\gamma}$. Consequently,
\begin{eqnarray*}
\Delta_{M+1}(t,s) &\leq& \Delta_{M+1}(t,\bar t)+ \Delta_{M+1}(\bar
s,s)+\Delta_{M}(\bar t,\bar s)\\
&&{}+|F_{M+1}(\bar s)-F_{M}(\bar s)|+|F_{M+1}(\bar t)-F_{M}(\bar t)| \\
&\le&
2(b-1)b^{-(M+1)\gamma}+2b \sum_{m=n+1}^{M} b^{-j\gamma}+2
b^{-(M+1)\gamma},
\end{eqnarray*}
so (\ref{eqn:10}) holds for $m=M+1$. Let $C_\gamma=2b/(1-b^{-\gamma
})$. Letting $M$ tend to infinity in~(\ref{eqn:10}) yields that $
\max_{k\geq1} |F_k(t)-F_k(s)| \leq C_\gamma b^{-(n+1)\gamma}$ for
all $n\geq n_1$ and $t,s\in T_*$ such that $|t-s|\leq b^{-n}$. Now, for
$t,s\in T_*$ with $|t-s|\leq b^{-n_1}$, there is a unique $n\geq n_1$
such that $b^{-(n+1)}\leq|t-s| <b^{-n}$ and $\max_{k\geq1}
|F_k(t)-F_k(s)| \leq C_\gamma b^{-(n+1)\gamma} \leq C_\gamma
|t-s|^\gamma$. The conclusion comes from the density of $T_*$ in
$[0,1]$ and the continuity of the $F_k$.
\end{pf*}
\noqed
\end{pf*}

%

\printaddresses

\end{document}